\documentclass[letterpaper, 12pt]{article}
\pdfoutput=1

\usepackage[nocompress]{cite}
\usepackage{jheppub}

\usepackage[table]{xcolor}

\usepackage{mathptmx} 

\usepackage{relsize} 
\usepackage{mathtools} 

\usepackage{graphicx}
\usepackage{epstopdf}
\usepackage{amsmath, amssymb}

\usepackage{bm}
\usepackage{caption}
\usepackage{subcaption}

\usepackage{multirow} 
\usepackage{longtable} 

\usepackage{romannum}

\usepackage[makeroom]{cancel} 

\usepackage{dsfont}

\usepackage{comment}

\usepackage{url}

\newcommand{\be}{\begin{eqnarray}}
\newcommand{\ee}{\end{eqnarray}}

\newcommand{\bn}{\begin{enumerate}}
\newcommand{\en}{\end{enumerate}}

\def\half{\frac{1}{2}}

\newcommand{\bea}{\begin{eqnarray}}
\newcommand{\eea}{\end{eqnarray}}

\def\half{\frac{1}{2}}

\allowdisplaybreaks

\title{Singlets in the tensor product of an arbitrary number of Adjoint representations of SU(3)}

\author[a, b]{Prarit Agarwal, }
\author[b, c]{June Nahmgoong}
\affiliation[a]{CRST and School of Physics and Astronomy, \\ Queen Mary University of London, London E1 4NS, UK}
\affiliation[b]{Department of Physics and Astronomy \& Center for Theoretical Physics\\ Seoul National University, Seoul 08826, Korea}
\affiliation[c]{School of Physics, Korea Institute of Advanced Study, Seoul 02455, Korea}

\abstract
{
We propose a set of 4 recurrence relations whose linear combination gives the number of group invariants, equivalently the dimension of the invariant subspace, in the tensor product of an arbitrary number of adjoint representations of the SU(3) Lie Group.

}


\preprint{QMUL-PH-20-02, KIAS-P20007}

\begin{document}
\maketitle


\section{Computing the dimension of invariant subspace}

Let $G$ be an arbitrary Lie group and $\{V_i, 1 \leq i \leq n\}$, be a set of $n$ arbitrary representations of $G$. We are interested in computing the dimension, $d_{inv}$, of the invariant subspace of the tensor product given by $V_1 \otimes \hdots \otimes V_n$. This can be computed by evaluating the integral
\be
\label{eq:SingletsIntegral}
d_{inv} = \oint_T d\mu_G \prod_{i = 1}^{n}\chi^{V_i}_G \ ,
\ee
where, $T \subseteq G$ is the maximal torus of $G$, $d\mu_G$ is the Haar measure and $\chi^{V_i}_G$ denotes the Weyl character of $V_i$. More explicitly, let $W$ be the Weyl group of G, $\Delta$ be the root system of $G$ and $\Delta^+$ be the set of positive roots, then
\be
\oint_T d\mu_G &=& \frac{1}{(2 \pi i)^r}\frac{1}{|W|} \oint_{|z_1| = 1} \hdots \oint_{|z_r| = 1} \frac{d z_1}{z_1} \hdots \frac{d z_r}{z_r} \prod_{\alpha \in \Delta}( 1 - z^{\alpha} ) \ , \\
z^{\alpha} &:=& \prod_{l=1}^{r} z_l^{\alpha_l} 
\ee 
where $\alpha_l$ denotes the $l$-th component of the root vector $\alpha$ (in the Dynkin basis) . Similarly, the character $\chi^{V_i}_G$ can be defined by 
\be
\chi^{V_i}_G &=& \sum_{\lambda \in V_i } z^{\lambda} \ , \\
z^{\lambda} &:=& \prod_{l=1}^{r} z_l^{\lambda_l} \ ,
\ee 
where $\lambda$ is a weight of $V_i$ and $\lambda_l$ denotes its $l$-th component. For example, the character of the fundamental irreducible representation of $SU(N)$ with the highest weight $[1, 0 , \hdots, 0]$  is given by:
\be
\chi^{[1, 0 , \hdots, 0]}_{SU(N)} = z_1 + \sum_{k=2}^{N-1} \frac{z_k}{z_{k-1}} + \frac{1}{z_{N-1}} \ .
\ee 
Similarly, the adjoint representation of $SU(N)$ has highest weight $[1,0, \hdots, 0, 1]$ with its character being given by 
 \be
 \chi^{[1, 0 , \hdots, 0, 1]}_{SU(N)} = (N-1) + \sum_{\alpha \in \Delta} z^\alpha \ .
 \ee 

Before bringing this section to a close, let us point out that from a computational point of view it turns out to be much more efficient to use a slightly different definition of the Haar measure \cite{Hanany:2008sb}, given by
\be
\label{eq:reducedHaarMeasure}
\oint_T d\mu_G &=& \frac{1}{(2 \pi i)^r} \oint_{|z_1| = 1} \hdots \oint_{|z_r| = 1} \frac{d z_1}{z_1} \hdots \frac{d z_r}{z_r} \prod_{\alpha \in \Delta^+}( 1 - z^{\alpha} ) \ .
\ee 
In the above definition, we restrict to the set of positive roots only. This removes the need to normalize by the order of the Weyl group in our previous definition. The measure given in \eqref{eq:reducedHaarMeasure} can be further simplified by observing that an application of the Weyl character formula to the trivial 1-dimensional representation implies that 
\be
\label{eq:HaarFurtherReduced}
\prod_{\alpha \in \Delta^+}( 1 - z^{\alpha} ) = \sum_{w \in W} { \rm sgn}(w) z^{w(\rho)-\rho} \ ,
\ee 
where $\rho$ is the Weyl vector defined by 
\be
\rho : = \half \sum_{\alpha \in \Delta^+} \alpha \ .
\ee
Simplifying \eqref{eq:reducedHaarMeasure} by substituting \eqref{eq:HaarFurtherReduced} leads to a further reduction in the number of computations that need to be performed to evaluate \eqref{eq:SingletsIntegral}.

\section{$n$-th tensor power of the adjoint representation of $SU(2)$}

Let us now apply the above method to compute the dimension, $d_{SU(2)}^{\rm inv} (n)$,  of the invariant subspace of the $n$-th tensor power of the adjoint representation of $SU(2)$. 

Following the discussion in the previous section, it is easy to see that
\be
d\mu_{SU(2)} &=& \frac{1}{2 \pi i}\frac{dz}{z}(1-z^2) \ , \\
\chi^{\rm adj.}_{SU(2)} &=& z^2 + 1 + z^{-2} \ , 
\ee  
which implies that
\be
d_{SU(2)}^{\rm inv}(n) &=& \frac{1}{2 \pi i}\oint_T d\mu_{SU(2)} \big(\chi^{\rm adj.}_{SU(2)} \big)^n \ , \\
&=& \frac{1}{2 \pi i} \oint_{|z|=1}\frac{dz}{z}(1-z^2)(1+z^2 + z^{-2})^n \ . \label{eq:SU2inv}
\ee 
Applying the residue theorem to the above integral we see that the dimension of the invariant subspace is given by 
\be
d_{SU(2)}^ {\rm inv}(n) &=& \text{the $z$-independent term in the expansion of } (1+z^2 + z^{-2})^n \nonumber \\
& & - \text{ the coefficient of $z^{-2}$ in the expansion of } (1+z^2 + z^{-2})^n \ .
\ee 
This implies
\be
d_{SU(2)}^{\rm inv}(n) = \mathlarger{\mathlarger{\sum}}_{r = 0,2, \hdots, n} \binom{n}{r}\binom{r}{r/2} - \mathlarger{\mathlarger{\sum}}_{r = 1,3, \hdots, n} \binom{n}{r}\binom{r}{(r-1)/2} \ .
\ee 
The first few explicit terms in this series are: $0, 1, 1, 3, 6, 15, 36, 91, 232, 603, 1585, 4213, \hdots$, which matches exactly with the sequence of the so called ``Motzkin sums'' i.e. \href{https://oeis.org/A005043}{A005043} on \href{https://oeis.org/}{oeis.org}. This observation has also been made in various forms in the comments on the afore mentioned online page, by David Callan, Andrew V. Sutherland and Samson Black. 

``Motzkin sums'' are known to obey the following recursion relation:
\be
a(n+1) = \frac{n}{n+2}\big(  2 a(n) + 3 a(n-1)\big) \ . \label{eq:SU2rec}
\ee
Let us now prove that the sequence of numbers given by the integral in \eqref{eq:SU2inv} indeed satisfies the recursion relation given by \eqref{eq:SU2rec}. In order to proceed let us define a generating function given as follows
\be
f_{SU(2)}(x) &=& \sum_{n=0}^{\infty}  d_{SU(2)}^{\rm inv}(n) x^n \, \label{eq:SU2GenSeries} \\
&=& \sum_{n=0}^{\infty} \frac{1}{2 \pi i} \oint_{|z|=1}\frac{dz}{z}(1-z^2)(1+z^2 + z^{-2})^n x^n \ , \\
&=& \frac{1}{2 \pi i} \oint_{|z|=1}\frac{dz}{z}\frac{1-z^2}{1-(1+z^2 + z^{-2}) x} \ ,  \label{eq:SU2gen1}
\ee 
where we will assume that $x$ is sufficiently small for the series to converge.  In order to evaluate $f_{SU(2)}(x)$ as defined in \eqref{eq:SU2gen1} it is somewhat easier to change the variable of integration from $z$ to $y = z^2$.  The generating function is then given by  
\be
f_{SU(2)}(x) =  \oint_{|y|=1}\frac{dy}{y}\frac{1-y}{1-(1+y + y^{-1}) x} \ .  \label{eq:SU2gen2}
\ee 
The integrand in the RHS of \eqref{eq:SU2gen2} has simple poles at 
\be
y &=& \frac{\sqrt{-3 x^2-2 x+1}-x+1}{2 x}\ , \ \text{and} \\
y &=&  -\frac{\sqrt{-3 x^2-2 x+1}+x-1}{2 x} \ .
\ee 
For $x<<1$, only the pole at $ -\frac{\sqrt{-3 x^2-2 x+1}+x-1}{2 x}$ lies within the unit circle and hence contributes to the contour integral in \eqref{eq:SU2gen2}. Upon evaluating the residue at this pole we therefore find that 
\be
f_{SU(2)}(x) = \frac{-1+3 x+\sqrt{1-2 x-3 x^2}}{2 x \sqrt{1-2 x-3 x^2}} \ . \label{eq:SU2GenClosedForm}
\ee 
Given the closed form expression \eqref{eq:SU2GenClosedForm} for $f_{SU(2)}(x)$ it is straightforward to show by direct substitution that $f_{SU(2)}(x)$ satisfies the following differential equation:
\be
(x-2x^2-3x^3)f'_{SU(2)}(x) + (1-3x^2) f_{SU(2)}(x) -1 = 0 \ , \label{eq:SU2GenDiff}
\ee 
where $f'_{SU(2)}(x)$ denotes the derivative of $f_{SU(2)}(x)$ with respect to x. Upon substituting the series expansion \eqref{eq:SU2GenSeries} of $f_{SU(2)}(x)$ in \eqref{eq:SU2GenDiff} we then arrive at the recursion relation in \eqref{eq:SU2rec}.

\section{$n$-th tensor power of the adjoint representation of $SU(3)$} \label{sec:SU3}
In this section we will compute the dimension, $d_{SU(3)}^{\rm inv}(n)$, of the invariant subspace of the $n$-th tensor power of the adjoint representation of $SU(3)$. In the process, we will give a prescription to decompose $d_{SU(3)}^{\rm inv}(n)$ into exactly 4 parts and obtain the recurrence relations satisfied by each of those parts. This will be the main result of this paper.  

Using \eqref{eq:reducedHaarMeasure} and \eqref{eq:HaarFurtherReduced}, we can write 
\be
d\mu_{SU(3)} &=& \frac{dz_1}{z_1} \frac{dz_2}{z_2}\left(1+z_1^3-\frac{z_1^2}{z_2}-\frac{z_2^2}{z_1}-z_1^2 z_2^2+z_2^3\right) \ .
\ee 
Also, the character of the adjoint representation of $SU(3)$ is  
\be
\chi^{\rm adj.}_{SU(3)} = 2+\frac{z_1}{z_2^2}+\frac{1}{z_1 z_2}+\frac{z_1^2}{z_2}+\frac{z_2}{z_1^2}+z_1 z_2+\frac{z_2^2}{z_1} \ . \label{eq:SU3AdjChar}
\ee
We therefore see that $d_{SU(3)}^{\rm inv}(n)$ can be obtained by evaluating the integral
\be
\label{eq:SU3SingletsIntegral}
d_{SU(3)}^{\rm inv}(n) &=& \oint_{|z_1| = 1}\oint_{|z_2| = 1} \frac{dz_1}{z_1} \frac{dz_2}{z_2} \left(1+z_1^3-\frac{z_1^2}{z_2}-\frac{z_2^2}{z_1}-z_1^2 z_2^2+z_2^3\right) \times \nonumber \\
& &\phantom{space1234}  \left(2+\frac{z_1}{z_2^2}+\frac{1}{z_1 z_2}+\frac{z_1^2}{z_2}+\frac{z_2}{z_1^2}+z_1 z_2+\frac{z_2^2}{z_1} \right)^n \ .
\ee

Notice, that the $\mathbb{Z}_2$ outer-automorphism symmetry of the $\mathfrak{su}(3)$ algebra implies that both $d\mu_{SU(3)}$ and $\chi^{\rm adj.}_{SU(3)} $, and therefore the integral in \eqref{eq:SU3SingletsIntegral}, are invariant under $z_1 \leftrightarrow z_2$. Using this exchange symmetry and applying the residue theorem to the integral in \eqref{eq:SU3SingletsIntegral}, we therefore surmise that $d_{SU(3)}^{\rm inv}(n)$ can written as
\be
d_{SU(3)}^{\rm inv}(n) = a_1(n) -  2 a_2(n) + 2 a_3(n) - a_4 (n) \ ,
\ee  
where, $a_1(n)$, $a_2(n)$, $a_3(n)$ and $a_4 (n)$ are defined as follows
\be
a_1(n) &:=& \text{the constant term in $\left(2+\frac{z_1}{z_2^2}+\frac{1}{z_1 z_2}+\frac{z_1^2}{z_2}+\frac{z_2}{z_1^2}+z_1 z_2+\frac{z_2^2}{z_1} \right)^n$}\ , \label{eq:a1} \\
a_2(n) &:=& \text{ the coefficient of $z_1 z_2^{-2}$  in $\left(2+\frac{z_1}{z_2^2}+\frac{1}{z_1 z_2}+\frac{z_1^2}{z_2}+\frac{z_2}{z_1^2}+z_1 z_2+\frac{z_2^2}{z_1} \right)^n$} \ , \label{eq:a2} \\
a_3(n) &:=& \text{ the coefficient of $z_1^{-3}$  in $\left(2+\frac{z_1}{z_2^2}+\frac{1}{z_1 z_2}+\frac{z_1^2}{z_2}+\frac{z_2}{z_1^2}+z_1 z_2+\frac{z_2^2}{z_1} \right)^n$}\ , \label{eq:a3} \\
a_4(n) &:=& \text{ the coefficient of $z_1^{-2} z_2^{-2}$  in $\left(2+\frac{z_1}{z_2^2}+\frac{1}{z_1 z_2}+\frac{z_1^2}{z_2}+\frac{z_2}{z_1^2}+z_1 z_2+\frac{z_2^2}{z_1} \right)^n$} \ . \label{eq:a4}
\ee 
While, it is straightforward to explicitly write $a_i (n), i = 1,2,3,4$ in terms of appropriate binomial coefficients, these expressions are rather complicated. Instead, we will suffice ourselves with listing the first few numbers appearing in the corresponding series and then proceed to obtain the recurrence relations describing them. 

\paragraph{Recursion relation for $a_1{(n)}$:} 
The first 10 values in this series are:1, 2, 10, 56, 346, 2252, 15184, 104960, 739162, 5280932, 38165260. Turns out this matches with the first 10 values in the sequence of Franel numbers \cite{Franel1, Franel2} i.e. the sequence \href{https://oeis.org/A000172}{A000172} on \href{https://oeis.org}{oeis.org}. As was shown by Franel, this sequence satisfies a recursion relation given as 
\be
\label{eq:a1rec}
a(n+1) = \frac{2+7n+7n^2}{(n+1)^2} a(n) + \frac{ 8 n^2}{(n+1)^2} a(n-1) \ . \label{eq:SU3a1rec}
\ee

\begin{figure}[t]
	\begin{center}
		\includegraphics[width=3.3in]{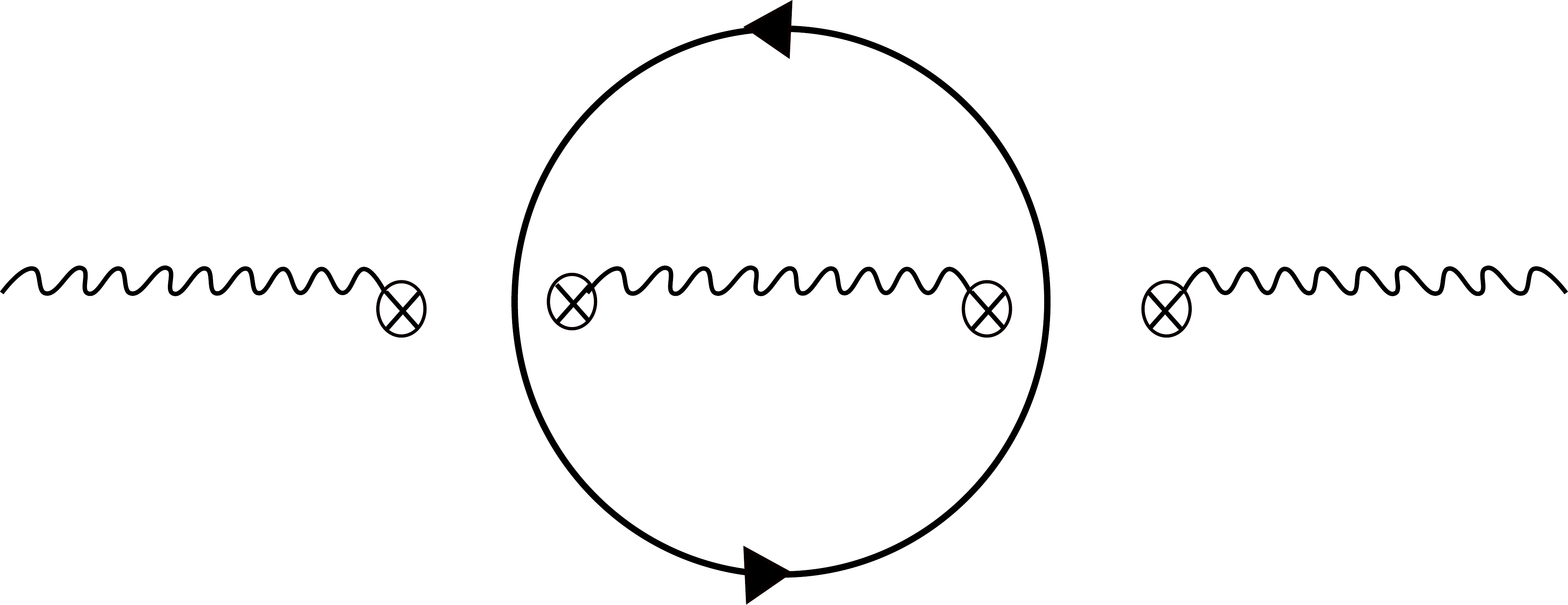}
		\caption{Contour for the integrals in section \ref{sec:SU3}. The branch-points are marked by $\otimes$ with the wavy lines being the branch cuts connecting them. }
		\label{fig:contour}
	\end{center}
\end{figure}

Let us now prove that $a_1(n)$ as defined in \eqref{eq:a1} indeed satisfies Franel's recursion relation as given in \eqref{eq:SU3a1rec}. As was the case in previous section, we will do so by defining a generating function for $a_1(n)$ and showing that it satisfies an appropriate differential equation.  Without further ado, let us define the following generating function:
\be
f_{SU(3),a_1}(x) &=& \sum_{n=0}^{\infty}  a_1(n) x^n \, \label{eq:SU3A1GenSeries} \\
&=& \sum_{n=0}^{\infty} \frac{1}{(2 \pi i)^2} \oint_{|z_2|=1} \oint_{|z_1|=1}\frac{dz_1}{z_1}\frac{dz_2}{z_2} (\chi^{\rm adj.}_{SU(3)})^n x^n \ , \\
&=& \frac{1}{(2 \pi i)^2} \oint_{|z_2|=1} \oint_{|z_1|=1}\frac{dz_1}{z_1}\frac{dz_2}{z_2}\frac{1}{1- \chi^{\rm adj.}_{SU(3)} x  } \ ,  \label{eq:SU3A1gen1}
\ee 
where $\chi^{\rm adj.}_{SU(3)}$ is defined in \eqref{eq:SU3AdjChar}. It is slightly simpler to analyze the above contour integrals by changing the integration variables to $y_1 = z_1^2/z_2, y_2 = z_2^2/z_1$. $f_{SU(3),a_1}(x)$ is then given by the following contour integral
\be
f_{SU(3),a_1}(x) = \frac{1}{(2 \pi i)^2} \oint_{|y_2|=1} \oint_{|y_1|=1}\frac{dy_1}{y_1}\frac{dy_2}{y_2} \frac{1}{1-(2 + y_1 + \frac{1}{y_1} + y_2 + \frac{1}{y_2} + y_1 y_2 + \frac{1}{y_1 y_2}) x  } \ . \nonumber \\
\ee 
Let us integrate with respect to $y_1$ first.  The simple poles in the $y_1$-plane are present at 
\be
y_1 &=& \frac{-x+y_2-2 x y_2-x y_2^2-\sqrt{x^2-2 x y_2+y_2^2-4 x y_2^2-2 x^2 y_2^2-2 x y_2^3+x^2 y_2^4}}{2 \left(x y_2+x y_2^2\right)} \ , \text{and} \nonumber \\ \\ 
y_1 &=&\frac{-x+y_2-2 x y_2-x y_2^2+\sqrt{x^2-2 x y_2+y_2^2-4 x y_2^2-2 x^2 y_2^2-2 x y_2^3+x^2 y_2^4}}{2 \left(x y_2+x y_2^2\right)} \ .
\ee  
For $x<<1$, only the pole at $\frac{-x+y_2-2 x y_2-x y_2^2-\sqrt{x^2-2 x y_2+y_2^2-4 x y_2^2-2 x^2 y_2^2-2 x y_2^3+x^2 y_2^4}}{2 \left(x y_2+x y_2^2\right)}$ lies withing the unit circle. Computing the residue at this point then gives
\be
f_{SU(3),a_1}(x) = \frac{1}{2 \pi i} \oint_{|y_2|=1}  \frac{dy_2}{\sqrt{x^2-2 x y_2+y_2^2-4 x y_2^2-2 x^2 y_2^2-2 x y_2^3+x^2 y_2^4}} \ . \label{eq:a1Period}
\ee 
The integrand above has 4 branch-points which in this case can be grouped into two pairs such that the product of the branch-points in any pair is 1. It therefore follows, the two of the branch points are inside the unit circle and two are outside. We will choose to connect the branch points inside the unit circle via a branch-cut. The other two branch points can be connected by another branch-cut that lies completely outside the unit circle. If one thinks of $y_1$ as the coordinate on $\mathbb{P}^1$ then the branch-points outside the unit circle can be connected by a branch-cut that also passes through the point at infinity. This is shown in figure \ref{fig:contour}. 
We therefore see that $f_{SU(3),a_1}(x)$ as given in \eqref{eq:a1Period} can be thought of as a period of the elliptic fibration given by 
\be
w^2 = x^2 y_2^4 -2 x y_2^3 +\left(-2 x^2-4 x+1\right) y_2^2-2 x y_2 +x^2 \ , \label{eq:SU3EllipticCurve}
\ee  
with $x$ being the coordinate on the base. It follows that $f_{SU(3),a_1}(x)$ satisfies the corresponding Picard-Fuchs differential equation. This is given by 
\be
x (x + 1) (8 x - 1) f_{SU(3),a_1}^{(2)}(x) + (24 x^2 + 14 x - 1) f_{SU(3),a_1}^{(1)}(x) + 
2 (4 x + 1) f_{SU(3),a_1}(x) = 0 \ , \nonumber \\ \label{eq:SU3A1diff}
\ee 
where $f_{SU(3),a_1}^{(n)}(x)$ denotes the $n$-th derivative of $f_{SU(3),a_1}(x)$ with respect to $x$. That this is indeed the case can be established in a straightforward manner by the usual brute-force approach, for e.g. as  described in \cite{PicardFuchsNotes}. Substituting the series expansion \eqref{eq:SU3A1GenSeries} in \eqref{eq:SU3A1diff} then gives the recursion relation \eqref{eq:a1rec}, thus proving that $a_1(n)$ is indeed given by Franel numbers.

\paragraph{Recursion relation for $a_2{(n)}$:}Let us now consider the first 10 values of $a_2(n)$ which are: 0, 1, 6, 39, 260, 1780, 12432, 88207, 633768, 4600566, 33680900. Unfortunately, there is no sequence in the \href{https://oeis.org}{oeis.org} database that matches this. However, inspired by the recurrence relation \eqref{eq:a1rec} satisfied by $a_1(n)$, we expect that $a_2(n)$ will also satisfy an analogous recurrence relation. A little bit of trial an error reveals that this indeed seems to be the case, with the recurrence relation we seek being given by
\be
a(n+1) =  \frac{(1+n) \left(8+30 n+49 n^2+21 n^3\right)}{n (2+n)^2 (1+3 n)} a(n) + \frac{8 n (1+n) (4+3 n)}{(2+n)^2 (1+3 n)} a(n-1) \ . \label{eq:a2rec}
\ee 
As before, we prove this by considering the generating function 
\be
f_{SU(3),a_2}(x) &=& \sum_{n=0}^{\infty}  a_2(n) x^n \, \label{eq:SU3A2GenSeries} \\
&=& \sum_{n=0}^{\infty} \frac{1}{(2 \pi i)^2} \oint_{|z_2|=1} \oint_{|z_1|=1}\frac{dz_1}{z_1}\frac{dz_2}{z_2} \frac{z_2^2}{z_1}(\chi^{\rm adj.}_{SU(3)})^n x^n \ , \\
&=& \frac{1}{(2 \pi i)^2} \oint_{|z_2|=1} \oint_{|z_1|=1}\frac{dz_1}{z_1}\frac{dz_2}{z_2}\frac{z_1^{-1} z_2^2}{1- \chi^{\rm adj.}_{SU(3)} x  } \ .  \label{eq:SU3A2gen1}
\ee 
Upon, changing the integration variables to $y_1 = z_1^2/z_2, y_2 = z_2^2/z_1$ and integrating with respect to $y_1$, we find that $f_{SU(3),a_2}(x)$ is given by the following contour integral:
\be
f_{SU(3),a_2}(x) = \frac{1}{2 \pi i} \oint_{|y_2|=1} dy_2  \frac{y_2}{\sqrt{x^2-2 x y_2+y_2^2-4 x y_2^2-2 x^2 y_2^2-2 x y_2^3+x^2 y_2^4}} \ . \label{eq:a2Period}
\ee 
Once again, the branch-cuts are as shown in figure \ref{fig:contour}. The recurrence relation in \eqref{eq:a2rec} can be shown to follow straightforwardly, if $f_{SU(3),a_2}(x)$ satisfies the differential equation:
\be
\left(-24 x^6-21 x^5+3 x^4\right)f_{SU(3),a_2}^{(4)}(x)+\left(-296 x^5-196 x^4+19 x^3\right) f_{SU(3),a_2}^{(3)}(x)+& \nonumber \\ \left(-968 x^4-436 x^3+19 x^2\right) f_{SU(3),a_2}^{(2)}(x)+\left(-848 x^3-208 x^2-2 x\right) f_{SU(3),a_2}^{(1)}(x)+& \nonumber \\ 
 \left(-112 x^2-8 x+2\right) f_{SU(3),a_2}(x) & = 0 \ . \nonumber \\ \label{eq:a2GenDiff}
\ee 
Substituting $f_{SU(3),a_2}(x)$ from \eqref{eq:a2Period} in the LHS of the above equation gives
\begin{align}
\text{L.H.S of \eqref{eq:a2GenDiff}}  = \frac{1}{2 \pi i} \oint_{|y_2|=1} dy_2 \frac{g_2(y_2, x)}{(x^2-2 x y_2+y_2^2-4 x y_2^2-2 x^2 y_2^2-2 x y_2^3+x^2 y_2^4)^{9/2}} \ , \nonumber \\ \label{eq:a2GenDiff2}
\end{align}
where $g_2(x)$ is a degree-16 polynomial in $y_2$ given by 
\begin{align}
g_2(y_2, x)=& (24 x^7-24 x^8) (y_2^{16} + y_2^2) + (-112 x^8-152 x^7-40 x^6) (y_2^{15}+y_2^3) + \nonumber \\
& (-264 x^8-1896 x^7+1056 x^6-90 x^5)(y_2^{14}+y_2^4) + \nonumber \\
& (-96 x^8-6384 x^7+2736 x^6-1344 x^5+282 x^4)(y_2^{13}+y_2^5) + \nonumber \\
& (936 x^8-9432 x^7+1344 x^6-658 x^5+484 x^4-268 x^3)(y_2^{12}+y_2^6) + \nonumber \\
&(1392 x^8-3048 x^7-5784 x^6+5520 x^5+1248 x^4-300 x^3+108 x^2)(y_2^{11}+y_2^7) + \nonumber \\
&(-648 x^8+11304 x^7-19296 x^6+8748 x^5- \nonumber \\ & \phantom{AAAAAAAAA} 324 x^4-1188 x^3+324 x^2-18 x)(y_2^{10} +y_2^8) + \nonumber \\ 
&(-2368 x^8+19168 x^7-27616 x^6+7648 x^5-\nonumber \\ & \phantom{AAAAAAAAA}2740 x^4-2312 x^3+320 x^2-44 x+2)y_2^{9} \ .
\end{align}
It is easy to show that the 1-form appearing in RHS of \eqref{eq:a2GenDiff2} is exact, being given by the following exterior derivative:
\begin{align}
 \frac{g_2(y_2, x)}{(x^2-2 x y_2+y_2^2-4 x y_2^2-2 x^2 y_2^2-2 x y_2^3+x^2 y_2^4)^{9/2}} = \nonumber \\ \partial_{y_2} \frac{h_2(y_2, x)}{(x^2-2 x y_2+y_2^2-4 x y_2^2-2 x^2 y_2^2-2 x y_2^3+x^2 y_2^4)^{7/2} (1+x)} \ ,
\end{align} 
where $h_2(y_2, x)$ is a degree-14 polynomial in $y_2$ being given by 
\begin{align}
h_2(x) = &  -2 \left(x^8+2 x^7\right) y_2^{14} + 2 x^5 \left(19 x^2+14 x-12\right) y_2^{13} + \nonumber \\
& 2 \left(7 x^8+56 x^7+43 x^6+6 x^5+40 x^4\right) y_2^{12} - \nonumber \\
& 2 \left(32 x^7-110 x^6-71 x^5+43 x^4+49 x^3\right) y_2^{11} - \nonumber \\
& 2 \left(21 x^8+170 x^7-349 x^6-202 x^5+211 x^4+49 x^3-27 x^2\right) y_2^{10} - \nonumber \\
&14 \left(3 x^7-84 x^6-30 x^5+17 x^4-16 x^3-13 x^2+x\right) y_2^9 + \nonumber \\
&2 \left(35 x^8+208 x^7+274 x^6-218 x^5+11 x^4+56 x^3+65 x^2-17 x+1\right) y_2^8 + \nonumber \\
& 2 \left(72 x^7-436 x^6-568 x^5+151 x^4-147 x^3+25 x^2-10 x\right)y_2^7 - \nonumber \\
& 2 \left(35 x^8+134 x^7+690 x^6+466 x^5-339 x^4+119 x^3-24 x^2\right)y_2^6 - \nonumber \\
& 2 \left(55 x^7+346 x^6+310 x^5-219 x^4+28 x^3\right) y_2^5 + \nonumber \\
& 2 \left(21 x^8+56 x^7+3 x^6-164 x^5+29 x^4\right)y_2^4 + \nonumber \\
& 2 \left(24 x^7+70 x^6-31 x^5\right) y_2^3 - \nonumber \\
& 14 x^6 \left(x^2+2 x-3\right)y_2^2 - \nonumber \\
& 14 x^7y_2^1 + \nonumber \\
& 2 x^8 \ .
\end{align}
Since the branch-cuts of $h_2(y_2, x)/(x^2-2 x y_2+y_2^2-4 x y_2^2-2 x^2 y_2^2-2 x y_2^3+x^2 y_2^4)^{7/2} (1+x)$ do not intersect the integration contour, it therefore follows that the RHS of \eqref{eq:a2GenDiff2} evaluates to zero by Stokes' theorem. We therefore see that the generating function $f_{SU(3),a_2}(x)$ satisfies the differential equation in \eqref{eq:a2GenDiff}. By substituting the series expansion given in \eqref{eq:SU3A2GenSeries}, it can  then be shown that $a_2(n)$ is given by the recursion relation in \eqref{eq:a2rec}.

\paragraph{Recursion relation for $a_3{(n)}$:}Similarly, the first 10 values of $a_3(n)$ are: 
0, 0, 2, 18, 144, 1100, 8280, 62034, 464576, 3484296, 26190900. As was the case for $a_2(n)$, we could not find any entry for a related sequence in the \href{https://oeis.org}{oeis.org} database. Nonetheless, empirically, we found that $a_3(n)$ satisfies the recursion relation given by
\be
a(n+1) = \frac{(1+n)^2 \left(-2+7 n+7 n^2\right)}{n \left(-6+n+4 n^2+n^3\right)} a(n)+\frac{8 n (1+n)^2}{-6+n+4 n^2+n^3} a(n-1) \ . \label{eq:a3rec}
\ee 
The proof that the above recursion relation is correct follows along the same lines as in the previous three cases i.e  consider the generating function for this sequence and prove that it satisfies an appropriate differential equation. The generating function is now given by 
\be
f_{SU(3),a_3}(x) &=& \sum_{n=0}^{\infty}  a_3(n) x^n \, \label{eq:SU3A3GenSeries} \\
&=& \sum_{n=0}^{\infty} \frac{1}{(2 \pi i)^2} \oint_{|z_2|=1} \oint_{|z_1|=1}\frac{dz_1}{z_1}\frac{dz_2}{z_2} z_1^3 (\chi^{\rm adj.}_{SU(3)})^n x^n \ , \\
&=& \frac{1}{(2 \pi i)^2} \oint_{|z_2|=1} \oint_{|z_1|=1}\frac{dz_1}{z_1}\frac{dz_2}{z_2}\frac{z_1^{3} }{1- \chi^{\rm adj.}_{SU(3)} x  } \ .  \label{eq:SU3A3gen1}
\ee
Changing the integration variables to $y_1 = z_1^2/z_2, y_2 = z_2^2/z_1$ and integrating with respect to $y_2$, we find that $f_{SU(3),a_3}(x)$ is given by the following contour integral:
\begin{align}
f_{SU(3),a_3}(x)=& -\frac{1}{2 \pi i} \oint_{|y_1|=1} dy_1 \times \nonumber \\ 
&\frac{y_1 \left(x-y_1+2 x y_1+x y_1^2+\sqrt{x^2-2 x y_1+y_1^2-4 x y_1^2-2 x^2 y_1^2-2 x y_1^3+x^2 y_1^4}\right)}{2 x \left(1+y_1\right) \sqrt{x^2-2 x y_1+y_1^2-4 x y_1^2-2 x^2 y_1^2-2 x y_1^3+x^2 y_1^4}} \ . \label{eq:a3Period}
\end{align}
The above integral has a simple pole at $y_1=-1$. We will therefore have to slightly deform our integration contour to either include or exclude this pole. While it is important to know  which way to deform the contour for purposes of computing the generating function exactly, our purpose here is to only establish the recursion relation for $a_3(n)$. As we will see shortly, this will not get affected by how we deform the contour.

For the recursion relation in \eqref{eq:a3rec} to be true, it is sufficient to show that the generating function $f_{SU(3),a_3}(x)$ satisfies the following  4-th order differential equation:
\begin{align}
\left(-8 x^6-7 x^5+x^4\right) f_{SU(3),a_3}^{(4)}(x)+\left(-96 x^5-63 x^4+6 x^3\right) f_{SU(3),a_3}^{(3)}(x)+& \nonumber\\ 
\left(-304 x^4-131 x^3+2 x^2\right) f_{SU(3),a_3}^{(2)}(x)+\left(-256 x^3-50 x^2-4 x\right) f_{SU(3),a_3}^{(1)}(x)+& \nonumber \\ 
\left(-32 x^2+2 x+4\right) f_{SU(3),a_3}(x)& = 0 \ .  \label{eq:a3GenDiff}
\end{align}
Substituting for $f_{SU(3),a_3}(x)$ from \eqref{eq:a3Period} in the LHS of \eqref{eq:a3GenDiff}, we find that 
\begin{align}
\text{L.H.S of \eqref{eq:a3GenDiff}}  = \frac{1}{2 \pi i} \oint_{|y_1|=1} dy_1 \frac{g_3(y_1, x)}{(x^2-2 x y_1+y_1^2-4 x y_1^2-2 x^2 y_1^2-2 x y_1^3+x^2 y_1^4)^{9/2}} \ , \nonumber \\ \label{eq:a3GenDiff2}
\end{align}
where $g_3(y_1,x)$ is an degree-16 polynomial in $y_1$ given by 
\begin{align}
g_3(y_1,x) =& 24 x^8 (y_1^{16} + y_1^3) + 24 x^2 \left(5 x^6+11 x^5-5 x^4\right) (y_1^{15}+y_1^4) + \nonumber \\ 
&  24 x^2 \left(6 x^6+45 x^5-36 x^4+15 x^3\right)(y_1^{14}+y_1^5) + \nonumber \\
& 24 x^2 \left(-10 x^6+97 x^5-43 x^4+34 x^3-15 x^2\right) (y_1^{13}+y_1^6) \nonumber \\
& 24 x^2 \left(-29 x^6+119 x^5+82 x^4+24 x^3-26 x^2+5 x\right) (y_1^{12}+y_1^7) \nonumber\\
& 24 x^2 \left(-9 x^6-18 x^5+207 x^4-90 x^3-36 x^2+18 x\right) (y_1^{11}+y_1^8) \nonumber \\
& 24 x^2 \left(36 x^6-254 x^5+179 x^4-287 x^3-83 x^2+17 x-2\right) (y_1^{10}+y_1^9) 
\end{align} 
One can check that the 1-form appearing in \eqref{eq:a3GenDiff2} is exact i.e. 
\begin{align}
\frac{g_3(y_1, x)}{(x^2-2 x y_1+y_1^2-4 x y_1^2-2 x^2 y_1^2-2 x y_1^3+x^2 y_1^4)^{9/2}} = \nonumber \\ \partial_{y_1} \frac{h_3(y_1, x)}{(x^2-2 x y_1+y_1^2-4 x y_1^2-2 x^2 y_1^2-2 x y_1^3+x^2 y_1^4)^{7/2}} \ , \label{eq:a3ExactForm}
\end{align} 
with $h_3(y_1, x)$ being a degree-14 polynomial in $y_1$:
\begin{align}
h_3(y_1,x) = &-6 x^6 y_1^{14} -6 \left(4 x^6-7 x^5\right) y_1^{13}-6 \left(3 x^6-2 x^5+11 x^4\right) y_1^{12} + \nonumber \\
& 6 \left(8 x^6-46 x^5-18 x^4+5 x^3\right) y_1^{11} + 42 \left(2 x^6-10 x^5+5 x^4+4 x^3\right) y_1^{10} - \nonumber \\
& 12 \left(8 x^5-48 x^4-19 x^3+2 x^2\right) y_1^9 - \nonumber \\
& 6 \left(14 x^6-38 x^5-59 x^4-20 x^3+8 x^2\right) y_1^8 - \nonumber \\
& 6 \left(8 x^6-46 x^5+10 x^4-13 x^3+4 x^2\right) y_1^7 + \nonumber \\
& 6 \left(3 x^6+30 x^5-19 x^4+8 x^3\right) y_1^6 + 6 \left(4 x^6+9 x^5-4 x^4\right) y_1^5 + 6 x^6 y_1^4 \ .
\end{align} 
We can therefore infer that since the contour of integration does not cross the branch-cuts of $h_3(y_1, x)/(x^2-2 x y_1+y_1^2-4 x y_1^2-2 x^2 y_1^2-2 x y_1^3+x^2 y_1^4)^{7/2}$, thus integration over the RHS of \eqref{eq:a3ExactForm} will evaluate to zero, thereby showing that $f_{SU(3),a_3}(x)$ satisfies the differential equation in \eqref{eq:a3GenDiff}. The recursion relation in \eqref{eq:a3rec} then follows by substituting the series expansion \eqref{eq:SU3A3GenSeries}. 

In order to tie up loose ends, recall that we had earlier claimed that for purposes of obtaining the recursion relation \eqref{eq:a3rec}, it does not matter how we deform the contour of integration in \eqref{eq:a3Period}. This is because the results so obtained differ from each other by the residue of the integrand appearing in \eqref{eq:a3Period}. One can check that this is $1/x$. Moreover, $1/x$ satisfies the differential equation in \eqref{eq:a3GenDiff}. Thus whether we include the pole at $y_1 = -1$ or not inside the contour of integration, either way the generating function will satisfy the same differential equation and therefore the coefficients of its series expansion will satisfy the same recursion relation. We see that the our choice of the deformed contour only changes the coefficient of the $x^{-1}$ term in the series expansion of the generating function which in any case is inconsequential, the coefficients of interest for us being coefficient of  $x^n, \ \forall n \geq 0$.

\paragraph{Recursion relation for $a_4{(n)}$:}Finally, we look at $a_4(n)$. The fist 10 numbers in this series are: 0, 0, 1, 12, 106, 860, 6735, 51912, 397180, 3029112, 23078100. These too don't seem to be present in the \href{https://oeis.org}{oeis.org} database. However, we found that the for all $0 \leq n \leq 100$, $a_4(n)$ satisfies the recurrence relation
\be
a(n+1) &=& \frac{(1+n) \left(36+116 n+91 n^2+21 n^3\right)}{(3+n)^2 \left(-6+n+3 n^2\right)} a(n) \nonumber \\ 
&&+ \frac{8 n^2 (1+n) \left(-2+7 n+3 n^2\right)}{(3+n)^2 \left(6-7 n-2 n^2+3 n^3\right)} a(n-1) \ .
\label{eq633_a4rec}
\ee 
As has been the case so far, in order to  prove this, it will be sufficient to show that the corresponding generating function satisfies an appropriate differential equation. The generating function for $a_4(n)$ is 
\begin{align}
f_{SU(3),a_4}(x)&=\sum_{n=0}^\infty a_4(n) x^n
\nonumber \\
&=\sum_{n=0}^\infty {1\over (2\pi i)^2}\oint_{|z_2|=1} \oint_{|z_1|=1} {dz_1\over z_1} {dz_2\over z_2} z_1^2 z_2^2 (\chi^\text{adj}_{SU(3)})^n x^n
\nonumber \\
&={1\over (2\pi i)^2}\oint_{|z_2|=1} \oint_{|z_1|=1} {dz_1\over z_1}{dz_2\over z_2}{z_1^2 z_2^2\over 1-\chi^\text{adj}_{SU(3)} x}.
\end{align}
Changing the integration variables to \(y_1=z_1^2/z_2\), \(y_2=z_2^2/z_1\) and integrating with respect to \(y_1\), we find that \(f_{SU(3),a_4}(x)\) is given by the following contour integral:
\begin{align}
f_{SU(3),a_4}(x)&={1\over 2\pi i }\oint_{|y_2|=1} dy_2
\Bigr[\frac{x y_2^2+(2 x-1) y_2+x}{2 x^2 \left(y_2+1\right){}^2}
+
\nonumber \\
&\frac{x^2 y_2^4+\left(2 x^2-4 x+1\right) y_2^2+x^2+2 (x-1) x y_2^3+2 (x-1) x
   y_2}{2 x^2 \left(y_2+1\right){}^2 \sqrt{x^2 y_2^4+\left(-2 x^2-4 x+1\right)
   y_2^2+x^2-2 x y_2^3-2 x y_2}}\Bigr]
   \label{eq652_fa4}
\end{align}
There is double pole at $y_2=-1$ with the corresponding residue being $-1/x^2$. Thus our choice of deforming the contour to include or exclude the pole will only change the coefficient of $x^{-2}$ term in series expansion of the generating function. Since the coefficient of $x^n\ ,  \forall n \geq 0$ are unaffected, the choice of deformed contour will not matter to us, as long as  $-1/x^2$ satisfies the same differential equation that governs the generating function. We, therefore, deform the contour of integration to exclude the pole at \(y_2=-1\).
\par 
For the recursion relation (\ref{eq633_a4rec}) to be true, it is sufficient to show that the generating function \(f_{SU(3),a_4}(x)\) satisfies the following 5-th order differential equation:
\begin{align}
(3x^5-21x^6-24x^7)f^{(5)}_{SU(3),a_4}(x)+
(31x^4-301x^5-440x^6)f^{(4)}_{SU(3),a_4}(x)&+
\nonumber \\
(55x^3-1166x^4-2400x^5)f^{(3)}_{SU(3),a_4}(x)+
(-38x^2-1182x^3-4384x^4) f^{(2)}_{SU(3),a_4}(x)&+
\nonumber \\
(22x-36x^2-2176x^3)f^{(1)}_{SU(3),a_4}(x)+
(32+36x-128x^2) f_{SU(3),a_4}(x)&=0.
\label{eq666_a4geneq}
\end{align}
Substituting \(f_{SU(3),a_4}(x)\) from (\ref{eq652_fa4}) in the LHS of (\ref{eq666_a4geneq}), we find that
\begin{align}
\text{L.H.S of }(\ref{eq666_a4geneq})={1\over 2\pi i }\oint_{|y_2|=1} dy_2
\Bigr[ {g_4(y_2,x)\over (x^2-2x y_2+y_2^2-4x y_2^2-2x^2 y_2^2-2x y_2^3+x^2 y_2^4)^{11/2} } -{6\over x}\Bigr],
\label{eq672_1form}
\end{align}
where \(g_4(y_2,x)\) is a degree-22 polynomial in \(y_2\) given by
\begin{align}
&g_4(y_2,x)=
-6 x^{10} (1+y_2^{22})
+66 x^9 (y_2+y_2^{21})
+(54 x^{10}+132 x^9-330 x^8) (y_2^2+y_2^{20})
\nonumber \\
&+(-396 x^9-1188 x^8+990 x^7) (y_2^{3}+y_2^{19})
+(-978 x^{10}+1956 x^9-1206 x^8+4752
   x^7-1980 x^6) (y_2^4+y_2^{18})
   \nonumber \\
&+(-1536 x^{10}+20826 x^9+18192
   x^8+4308 x^7-11088 x^6+2772 x^5)(y_2^5+y_2^{17})
\nonumber \\
&+(4290 x^{10}+45072 x^9+71940
   x^8-54024 x^7-10596 x^6+16632
   x^5-2772 x^4)(y_2^{6}+y_2^{16})
\nonumber \\
&+(9216 x^{10}-6672 x^9+114864
   x^8-178176 x^7+27336 x^6+28152
   x^5-16632 x^4+1980 x^3)(y_2^{7}+y_2^{15})
\nonumber \\
&+(-7452 x^{10}-148272 x^9+108060
   x^8-150984 x^7+108228 x^6
   \nonumber \\
   &\quad +32304
   x^5-40584 x^4+11088 x^3-990 x^2)(y_2^8+y_2^{14})
\nonumber \\
&+(-23040 x^{10}-156060 x^9+22512
   x^8+102828 x^7+122064 x^6
   \nonumber \\
   &\quad +29772
   x^5-58632 x^4+28716 x^3-4752
   x^2+330 x)(y_2^9+y_2^{13})
\nonumber \\
&+(4092 x^{10}+101112 x^9-178464
   x^8+255552 x^7-85668 x^6
   \nonumber \\
   &\quad -37224
   x^5-79332 x^4+44352 x^3-10098
   x^2+1188 x-66)(y_2^{10}+y_2^{12})
\nonumber \\
&+(30720 x^{10}+284472 x^9-308760
   x^8+250692 x^7-256656 x^6
   \nonumber \\
   &\quad -97968
   x^5-94848 x^4+49488 x^3-12672
   x^2+1716 x-132+\frac{6}{x})y_2^{11}.
\end{align}
One can check that the 1-form in (\ref{eq672_1form}) is exact i.e.
\begin{align}
&{g_4(y_2,x)\over (x^2-2x y_2+y_2^2-4x y_2^2-2x^2 y_2^2-2x y_2^3+x^2 y_2^4)^{11/2} }=
\nonumber \\
\partial_{y_2} & {h_4(y_2,x)\over{ (1+x)^2 (x^2-2x y_2+y_2^2-4x y_2^2-2x^2 y_2^2-2x y_2^3+x^2 y_2^4)^{9/2} } }
\end{align}
where \(h_4(y_2,x)\) being a degree-19 polynomial in \(y_2\) given by:
\begin{align}
&h_4(y_2,x)= 
(8 x^{11}-16 x^{10}-8 x^9)+
(-78 x^{10}+132 x^9+66 x^8)y_2
\nonumber \\
&+(-72 x^{11}+702 x^9-324 x^8-234 x^7) y_2^2
+(482 x^{10}+352 x^9-3070 x^8+12
   x^7+456 x^6)y_2^3
\nonumber \\
&+(288 x^{11}+432 x^{10}-2592 x^9-5364
   x^8+5760 x^7+1260 x^6-504 x^5)y_2^4
\nonumber \\
&+(-1344 x^{10}-2796 x^9+4116 x^8+18552
   x^7-3648 x^6-2268 x^5+252 x^4)y_2^5
\nonumber \\
&+(-672 x^{11}-1936 x^{10}+6040
   x^9+25096 x^8+27180 x^7-20100
   x^6-3280 x^5+1428 x^4+84 x^3)y_2^6
\nonumber \\
&+(2112 x^{10}+12684 x^9+27732
   x^8+12408 x^7-45564 x^6+1128
   x^5+5616 x^4+324 x^3-216 x^2)y_2^7
\nonumber \\
&+(1008 x^{11}+4368 x^{10}-768
   x^9-10284 x^8-31488 x^7-72516
   x^6
   \nonumber \\
   &\quad +9276 x^5+10212 x^4-852
   x^3-972 x^2+144 x)y_2^8
\nonumber \\
&+(-1804 x^{10}-23060 x^9-49620
   x^8-59612 x^7-69196 x^6
   \nonumber \\
   &\quad +27444
   x^5+14116 x^4-3520 x^3-1902
   x^2+520 x-46)y_2^9
\nonumber \\
&+(-1008 x^{11}-5904 x^{10}-18828
   x^9-35592 x^8-38580 x^7-32808
   x^6
   \nonumber \\
   &\quad +54696 x^5+24240 x^4-3360
   x^3-2628 x^2+798 x-96+\frac{6}{x})y_2^{10}
\nonumber \\
&+(468 x^{10}+12348 x^9+12204 x^8-12924
   x^7-30204 x^6
   \nonumber \\
   &\quad +45108 x^5+29268
   x^4+576 x^3-2790 x^2+648 x-54)y_2^{11}
\nonumber \\
&+(672 x^{11}+4880 x^{10}+25024
   x^9+34628 x^8+608 x^7-50420
   x^6
   \nonumber \\
   &\quad -1748 x^5+16180 x^4+3524
   x^3-1836 x^2+216 x)y_2^{12}
\nonumber \\
&+(528 x^{10}+5340 x^9+18084 x^8+19464
   x^7-26556 x^6-20616 x^5+1296
   x^4+2772 x^3-504 x^2)y_2^{13}
\nonumber \\
&+(-288 x^{11}-2352 x^{10}-9960
   x^9-5352 x^8+19548 x^7+6012
   x^6-7392 x^5-2268 x^4+756 x^3)y_2^{14}
\nonumber \\
&+(-496 x^{10}-5660 x^9-9788 x^8+2968
   x^7+6608 x^6+756 x^5-756 x^4)y_2^{15}
\nonumber \\
&+(72 x^{11}+576 x^{10}+360 x^9-2916
   x^8-2304 x^7+252 x^6+504 x^5)y_2^{16}
\nonumber \\
&+(138 x^{10}+672 x^9+282 x^8-324
   x^7-216 x^6)y_2^{17}
\nonumber \\
&+(-8 x^{11}-48 x^{10}+30 x^9+108
   x^8+54 x^7)y_2^{18}
   +(-6 x^{10}-12 x^9-6 x^8)y_2^{19}.
\end{align}
Since the second term of (\ref{eq672_1form}) given by \(-{6/ x}\), does not contribute to the integral, we can conclude that the RHS of (\ref{eq672_1form}) is zero. Then the generating function \(f_{SU(3),a_4}(x)\) satisfies the differential equation in (\ref{eq666_a4geneq}), thereby proving the recursion relation \eqref{eq633_a4rec} for \(a_4(n)\) . 

Coming back to the question of our choice of the deformed contour, one can easily check that $-1/x^2$ .i.e. the residue of integrand in the RHS of \eqref{eq652_fa4} satisfies the differential equation \eqref{eq666_a4geneq} thereby confirming that the choice of the deformed contour will be inconsequential for our purposes. 

\section{Discussion: $n$-th tensor power of the adjoint representation of $SU(N)$}

At this point it is natural to wonder if similar recursion relations can also also be established for the dimension of the invariant subspace in the $n$-th tensor power of the adjoint representation of $SU(N)$, for general $N$. A more conservative goal would be to establish an analogous system of recursion relations for $N=4$. This is work in progress. 

Alternatively, the dimension of the invariant subspace in the $n$-th power of the $SU(N)$ adjoint representation is given by the number of distinct ways to distribute $n$ matrices $T_i, 1 \leq i \leq n$ into different groups and take traces over their products with the condition ${\rm T_i} = 0 \ , \forall i$, which accounts for the fact that $T_i$ is an arbitrary generator of $SU(N)$. Additionally, one also has to take into consideration that ${\rm Tr}T_{i_1} \hdots T_{i_r} \ , r > N$ is not an independent tensor and can be written as a linear combination of   ${\rm Tr}T_{i_1} \hdots T_{i_s} \ , \forall s \leq N$. These relations, proliferate as $n$ and $N$ become large thereby making the counting-process a highly computationally intensive task. As a corollary, it follows that for $n \leq N$, the invariant subspace in the $n$-th power of the adjoint representation is independent of $N$ and is known to be equal to the number of derangements of a set of size $n$ \cite{Cvitanovic:1976am}. This is can also be confirmed through explicit computations for some low values of $n$ and $N$. The corresponding results are shown in table \ref{tab:SUNAdjProd}.
\begin{table}[t!]
	\centering
	\begin{tabular}{|c|c|c|c|c|c|c|c|c|}\hline
		\text{n} & SU(2) & SU(3)& SU(4)& SU(5)& SU(6)& SU(7) & SU(\(\infty\))\\ \hline 
		2 & \cellcolor{blue!25}1 & \cellcolor{blue!25}1 & \cellcolor{blue!25}1 & \cellcolor{blue!25}1 & \cellcolor{blue!25}1 & \cellcolor{blue!25}1 & \cellcolor{blue!25}1 
		\\ \hline
		3 & 1 & \cellcolor{blue!25}2 & \cellcolor{blue!25}2 & \cellcolor{blue!25}2 & \cellcolor{blue!25}2 & \cellcolor{blue!25}2 & \cellcolor{blue!25}2 
		\\ \hline
		4 & 3 & 8 & \cellcolor{blue!25}9 & \cellcolor{blue!25}9 & \cellcolor{blue!25}9 & \cellcolor{blue!25}9 & \cellcolor{blue!25}9 
		\\ \hline
		5 & 6 & 32 & 43 & \cellcolor{blue!25}44 & \cellcolor{blue!25}44 & \cellcolor{blue!25}44 & \cellcolor{blue!25}44
		\\ \hline
		6& 15 & 145 & 245 & 264 & \cellcolor{blue!25}265 & \cellcolor{blue!25}265 & \cellcolor{blue!25}265
		\\ \hline
		7 & 36 & 702 & 1557 & 1824 & 1853 & \cellcolor{blue!25}1854 & \cellcolor{blue!25}1854 
		\\ \hline
		8 & 91 & 3598 & 10829 & 14210 & 14791 & 14832 & \cellcolor{blue!25}14833 
		\\ \hline
	\end{tabular}
	\caption{Dimension of the invariant subspace in the $n$-th tensor power of the adjoint representation of $SU(N)$ Lie groups. Notice that the dimensions given in the shaded part stay constant along the rows.}
	\label{tab:SUNAdjProd}
\end{table}

\acknowledgments
We thank Sunjin Choi, Seok Kim, Costis Papageorgakis and Matthieu Sarkis for useful and enlightening discussions. The work of P.A. is supported by the Royal Society through a Research Fellows Enhancement Award, grant no. RGF\textbackslash EA\textbackslash 181049. The work of J.N. is supported by a KIAS Individual Grant PG076401.

\bibliographystyle{jhep}
\bibliography{ADN1}

\end{document}